\documentclass[11pt]{article}
\usepackage[english]{babel}
\usepackage{amsmath}
\usepackage{amsthm}
\usepackage{amssymb}
\newtheorem{theorem}{Theorem}[section]

\newtheorem{definition}[theorem]{Definition}
\newtheorem{algo}[theorem]{Algorithm}
\newtheorem{algorithm}[theorem]{Algorithm}
\newtheorem{lemma}[theorem]{Lemma}
\newtheorem{prop}[theorem]{Proposition}

\newtheorem{example}[theorem]{Example}
\newtheorem{corollary}[theorem]{Corollary}

\newtheorem{remark}[theorem]{Remark}

\newcommand{\pp}{\partial}

\newcommand{\cc}{{\cal C}}
\newcommand{\G}{{\cal G}}

\newcommand{\ccc}{{\mathbb{C}}}
\newcommand{\nn}{{\mathbb{N}}}
\newcommand{\rr}{{\mathbb{R}}}

\def\TTo#1{\mathop{\longrightarrow}\limits^{#1}}
\def\TToupdown#1#2{\mathop{\rightleftarrows}\limits^{#1}_{#2}}
\def\V{{\cal V}}
\def\TT{{\mathbb T}}
\def\NF{{\rm NF}_{\sigma}\ }

\def\GB{{{\rm  Gr\"obner Basis\ }}}
\def\GGBB{{{\rm Gr\"obner Bases\ }}}
\def\eGB{{{\rm  Gr\"obner Basis}}}
\def\eGGBB{{{ Gr\"obner Bases}}}
\def\terms#1{{\mathbb{T}^{#1}}}
\def\LT{{\rm LT}}
\def\cocoa{\mbox{\rm
   C\kern-.13em o\kern-.07 em C\kern-.13em o\kern-.15em A\ }}
   \def\ecocoa{\mbox{\rm
   C\kern-.13em o\kern-.07 em C\kern-.13em o\kern-.15em A}}
\def\coala{\mbox{\rm
   C\kern-.13em o\kern-.07 em A\kern-.13em l\kern-.15em A\ }}
\def\ecoala{\mbox{\rm
   C\kern-.13em o\kern-.07 em A\kern-.13em l\kern-.15em A}}

\title{\bf Computational Methods for the Construction of a Class of Noetherian Operators}
\author{
Alberto Damiano\\
Department of Mathematics\\
Charles University\\
Sokolovsk\'a 83, Praha, Czech Republic\\
alberto@tlc185.com\and
Irene Sabadini\\
Dipartimento di Matematica\\
Politecnico di Milano\\
Via Bonardi, 9\\
20133 Milano, Italy\\
sabadini@mate.polimi.it\and
Daniele C. Struppa\\
Department of Mathematics and Computer Sciences\\
Chapman University \\
Orange, CA 92866 \\
struppa@chapman.edu}
\date{  }
\begin{document}
\maketitle

\begin{abstract}
\noindent This paper presents some algorithmic techniques to
 compute explicitly the noetherian operators associated to a class of ideals
and modules over a polynomial ring. The procedures we include in
this work can be easily encoded in computer algebra packages such as
CoCoA \cite{cocoa}.
% We also
%consider some examples in positive dimension and we show how to
%apply the same techniques, with some modifications, to construct
%the operators.
\end{abstract}

\section{Introduction}
The Ehrenpreis--Palamodov Fundamental Principle, \cite{ehrenpreis} and \cite{palamodov},
states the following:
\begin{theorem}\label{t3.1}
Let $p_1(D),\ldots ,p_r(D)$ be linear constant coefficients
partial differential operators in $n$ variables. Then there are
algebraic varieties $V_1,\ldots ,V_t$ in $\ccc^n$ and differential
operators $\pp_1,\ldots ,\pp_t$ with polynomial coefficients, such that every function
$f\in\cc^{\infty}(\rr^n)$ satisfying
 $$
p_1(D)f=\ldots =p_r(D)f=0 $$ can be represented as
\begin{equation}\label{3.3}
f(x)=\sum_{j=1}^t\int_{V_j}\pp_j(e^{ix\cdot z})d\nu_j(z),
\end{equation}
for suitable Radon measures $d\nu_j$.
\end{theorem}
The collection
$${V}=\{(V_1,\pp_1);(V_2,\pp_2);\ldots;(V_t,\pp_t)\}$$
is said to be a multiplicity variety and Theorem \ref{t3.1} is equivalent to the following
strengthening of the classic Nullstellensatz:
\begin{theorem}\label{null1} Let $I$ be an ideal of $R$. There exists a multiplicity
 variety $V$ such that a polynomial $f$
belongs to $I$ if and only if $\pp_jf_{|V_j}=0$ for every
$j=1,\ldots,t$.
\end{theorem}
The operators $\pp_1,\ldots ,\pp_t$ are called, in Palamodov's
terminology, noetherian operators because their construction
relies essentially on a theorem of M. Noether on a membership
criterium for polynomial submodules (see e.g. \cite{palamodov}
pp.161, 162). The nature of the original proof of the Fundamental
Principle is essentially existential and therefore the question of
the explicit construction of such operators is of great interest
whenever we consider a concrete application of the Fundamental
Principle.
Note that if $I$ is the ideal generated by the polynomials $p_1,\ldots, p_r$
and if
$$I=\it{Q}_1\cap\dots\cap\it{Q}_t$$
is its primary decomposition, then the varieties $V_j$ which appear
in theorem \ref{t3.1} are simply given by the algebraic sets ${\cal
V}(Q_j)$. The information on the multiplicity of each of them is
left to the operators $\pp_j$.
%We point out that

In this paper we build on some recent results in the construction of
noetherian operators \cite{hop,mmm,oberst,oberst1,sturmfels} and
provide some new algorithms which allow the automatic construction
of these operators at least in some rather large class of cases. We
include several experiments using algorithms implemented on \cocoa
\cite{cocoa}.

In section 2 we quickly review the fundamental tools from
computational algebra (mostly the theory of \eGGBB).  The core of
the paper is section 3 where we deal with case of zerodimensional
ideals and where we present several explicit algorithms.  A final
section deals with the case of ideals of positive dimension.

Executable versions of the algorithms discussed in this paper have
been explicitly written for \cocoa and are freely available at
\begin{verbatim} http://www.tlc185.com/coala \end{verbatim}

The authors are grateful to F. Colombo, S. Hosten, and B. Sturmfels
for their many useful suggestions and comments. The first author is
grateful to George Mason University and to the Eduard \v{C}ech
Center for the financial support.

\section{Computational Algebra Tools}
Throughout this paper, we will work in the ring
$R=\ccc[x_1,\dots,x_n]$ of polynomials in $n$ variables with
complex coefficients; we will think of $R$ as the ring of symbols
for the differential operators we are studying. Even though we
consider differential operators with constant coefficients, the
Fundamental Principle shows that noetherian operators have, in
general, polynomial coefficients; we will use the symbol $A_n$ to
denote the Weyl Algebra
$\ccc[x_1,\dots,x_n,\pp x_1,\dots,\pp x_n]$ of such operators. Here, and throughout the paper, the symbol
$\pp x$ will be used as a shortcut for $\frac{\pp}{\pp x}$.\\
\noindent Using the notation introduced in \cite{kr}, we will
denote the monoid of power products in $R$ by $\TT^{n}$ and the module monoid of power products in $R^s$ by
$${\TT^{n}}\langle e_1,\dots,e_s\rangle=\{te_i\ |\ t\in\TT^{n},\, i=1,\dots,s\}$$ where
$e_i$ is the $i$-th element of the canonical basis of $R^s$. All the definitions to follow are given for ideals
but can be extended in straightforward fashion to the case of modules \cite{kr}.
A term ordering $\sigma$ on $\TT^{n}$ is a total ordering on power products with the following two properties:\\
\par\noindent
I)$\;\,$ if $t_1>_{\sigma}t_2$ and $t \in \terms{n}$ then $t\cdot t_1>_{\sigma}t\cdot t_2$;\\
II) if $t\in\terms{n}$ and $s\in\terms{n}$ then $s\cdot t>_{\sigma}t$.\\
\par\noindent
The {\it leading term ideal} associated to $I$ with respect to
$\sigma$ is the ideal generated by all the leading terms of elements
of $I$, and will be indicated by
$${\rm LT}_\sigma(I)=(\{ {\rm LT}_\sigma(f)|f\in I\}).$$ More in
general, the {\it leading term ideal} associated to a subset $G$ of
$R$ will be written as ${\rm LT}_\sigma(G)=(\{ {\rm
LT}_\sigma(f)|f\in G\})$. Note that ${\rm LT}_\sigma(G)={\rm
LT}_\sigma(I)$ if and only if the set $G$ is a \GB for the ideal
$I$, this being the main
characterization of a \eGB.\\
The algorithm which associates to an ideal $I$ of $R$  its
\GB ${\cal G}_{\sigma}(I)$  is the core algorithm of the theory of
\GGBB and can be found for example in \cite{kr}, theorem 2.5.5.
Another key tool in computational algebra is the division
algorithm (see again \cite{kr}, theorem 1.6.4) which can be
performed to generate the remainder  of a polynomial
 with respect to a set of generators of $I$.
Note that the remainder of a polynomial depends on the set of
generators chosen for $I$ (in fact, it even depends on their
order). The fundamental property of
\GGBB is that such a remainder is zero if and only if the
polynomial belongs to the ideal. For this reason the remainder
calculated with respect to a \GB is called the {\it normal form}
of a polynomial. \\
Given a polynomial $f\in I$ and a
term ordering $\sigma$, we will denote by $\NF(f)$ the normal
form of $f$ with respect to the $\sigma$--\GB of $I$ (the same
notation is used for modules). An equivalent way to compute a
remainder is using {\it rewrite rules} (see \cite{kr} section
2.2). Given a polynomial $g\in R$, we say that a polynomial $f_1$
rewrites to $f_2$ with respect to the rewrite rule $\TTo{{g}}$
(and this is indicated by $f_1\TTo{{g}}f_2$) if there exists a
monomial $m$ in $R$ such that $f_2=f_1-mg$ and ${\rm
LT}_\sigma(mg)$ is not in the support of $f_2$. This is also
called a one-step reduction. We can rewrite a polynomial using a
set of elements ${\cal G}=\{g_1,\dots,g_s\}$ by performing a
one-step reduction with each of the $g_i$'s, in that order. We
will denote by $\TTo{{\cal G}}$ the transitive closure of the
relations $\TTo{{g_1}},\dots,\TTo{{g_s}}$. This relation is called
{\it rewrite relation} or {\it rewrite rule}. By applying a
sequence of one-step reductions to a polynomial $f$ using the
elements in ${\cal G}$ we then obtain a remainder of $f$ with
respect to $\{g_1,\dots,g_s\}$. In particular if ${\cal G}$ is a
\GB we have that $f$ rewrites to its normal form, i.e.
$f\TTo{{\cal G}}\NF(f)$.

We now introduce
some definitions about elimination theory and term orderings (see \cite{kr},
section 3.4, for details on this topic) which will be necessary in the last section.

\begin{definition}\label{elimination}
Let $R=\ccc[x,t]$ where $x=(x_1,\dots,x_{n-d})$, $t=(t_1,\dots,t_d)$. A term ordering $\sigma$
on $\TT^n$ is called an {\rm elimination ordering} with respect to $x$ if every element $f\in
R$ whose leding term is contained in $\ccc[t]$ is such that $f\in\ccc[t]$. In other words,
$$\forall f\in R,\quad\LT_\sigma(f)\in\ccc[t]\,\Rightarrow\, f\in \ccc[t].$$
\end{definition}
The reason why such a term ordering is called an {\it elimination}
ordering is that it allows to eliminate the variables $x$ from an
ideal, i.e. it allows to compute $I\cap\ccc[t]$. To do this, it
suffices to compute a \GB with respect to any elimination ordering
as in definition \ref{elimination} and then keep only the elements
that do not contain any monomials in $x$. Such elements actually
form a \GB for the ideal $I\cap\ccc[t]$. It can be easily checked
that ${\tt Lex}$, the lexicographic term ordering on $\terms{n}$, is
an elimination ordering with respect to any "initial" subset of
variables, i.e. with respect to any subset of the type
$\{x_1,\dots,x_k\}$ in $\ccc[x_1,\dots,x_n]$, with $k\leq n$. A
class of term orderings that satisfy the elimination property and
that we are going to use for our goal of computing the noetherian
operators in $\ccc(t)[x]$ are the so called {\it product orderings}.
\begin{definition}\label{productordering}
Let $R=\ccc[x,t]$ as before and let $\sigma_x$ and $\sigma_t$ be two
term orderings on the set of terms $\TT_x=\{x^a\;|\;a\in\nn^{n-d}\}$
and $\TT_t=\{t^b\;|\;b\in\nn^d\}$ respectively. The {\rm product
ordering} $\sigma_x\cdot \sigma_t$ is defined by
$$x^at^b>_{\sigma_x\cdot \sigma_t}x^ct^d\,\Leftrightarrow\, x^a>_{\sigma_x}x^c\, {\it or}\,
(x^a=x^c\,{\it and}\,t^b>_{\sigma_t}t^d).$$
\end{definition}
It is immediate to show that the product ordering defined above is
an elimination ordering with respect to $x$, no matter what the
choice of $\sigma_x$ and $\sigma_t$ is. Elimination orderings are
usually slow when it comes to \GB computations, in particularly
${\tt Lex}$ is known to be one of the slowest. Product orderings are
then introduced to perform better. One can in fact define a "fast"
term ordering (such as ${\tt DegRevLex}$) on each of the two subsets
of variables, and then take the product. The following lemma will be
useful later in the paper.

\begin{lemma}\label{lemmaproduct}
Let $R=\ccc[x,t]$ be a polynomial ring equipped with a product ordering $\sigma$ of the type
$\sigma_x\cdot \sigma_t$ as in definition \ref{productordering}. Let $I$ be an ideal of $R$ and let
${\cal G}=(g_1,\dots,g_s)$ be a $\sigma$--\GB for $I$. Consider the extended ideal $IR_d$ in $R_d=\ccc(t)[x]$
endowed with the term ordering $\sigma_x$.
Then ${\cal G}$ forms a \GB for $IR_d$ with respect to $\sigma_x$.
\end{lemma}
\begin{proof}
Denote with $x^{a_i}t^{c_i}$ the leading term of $g_i$, where $a_i\in\nn^{n-d}$ and $c_i\in\nn^d$, $i=1,\dots,s$. From the fact that we chose a product ordering $\sigma$, it
follows  that once we view $g_i$ as an element of $IR_d$,
its leading term is $x^{a_i}$. In other words, $\LT_{\sigma_x}(g_i)=x^{a_i}$ in $R_d$. Consider a polynomial $f$
in $IR_d$. The set ${\cal G}$ still forms a set of generators for the extended ideal, so $f$ can be written as an
$R_d$-linear combination of the $g_i$'s. Moreover, supposing $f$ monic, we can write $f$ as
$$f=x^a+\sum_{b}p_b(t)x^b,\;{\rm where}\, b\in \nn^{n-d}\,{\rm and}\,x^a>_{\sigma_x}x^b\,\forall\,b.$$
Consider the product $D(t)$ of all the denominators of the coefficients $p_b(t)$ in $f$. Then $D(t)f$ is a
polynomial in $R$ and it is still a combination of the elements of ${\cal G}$, so $D(t)f\in I$. Because of the fact that $\sigma$ is a product order, the leading term of $D(t)f$ is simply the leading term of $f$ multiplied by some power of $t$, i.e. $\LT_\sigma(D(t)f)=x^at^c$ for some $c\in\nn^d$. Hence, ${\cal G}$ being a \GB for $I$, $x^at^c$ is a multiple of one of the leading terms of its elements, say $x^{a_1}t^{c_1}$ modulo a change on the order in ${\cal G}$. This means that there exist $\alpha\in\nn^{n-d}$ and $\gamma\in\nn^d$ such that
$$x^at^c=x^{\alpha}t^{\gamma}x^{a_1}t^{c_1}$$
which means that $x^a$ is a multiple of $x^{a_1}$, and this concludes the proof.
\end{proof}

\section{The Zerodimensional case}\label{zerodimsection}

In this section, $I$ is a primary zerodimensional ideal, i.e. the
algebraic set ${\cal V}(I)$ is a finite union of points in $\ccc^n$.
Since a zerodimensional primary ideal is associated to a single
point of the variety ${\cal V}(I)$ we can always assume, with a
change of coordinates, that ${\cal V}(I)=\{(0,\dots,0)\}$, or
equivalently that $\sqrt{I}=(x_1,\dots,x_n)$.

\subsection{Closed Differential Conditions}\label{closure}

A first complete description of the differential condition
characterizing a zerodimensional primary ideal centered in zero
has been done in \cite{mmm}: we briefly recall the main
notations and definitions of that paper. We will denote with
$D(i_1,\dots,i_n):R\to R$ the differential operator defined by:
$$D(i_1,\dots,i_n)=\frac{1}{i_1!\cdots i_n!}{\pp x_1^{i_1}\cdots\pp
x_n^{i_n}}, \quad i_j\in\nn\; ,\,{\rm for}\,{\rm all}\, j=1,\dots,n,$$ or, alternatively, if
$t=x_1^{i_1}\cdots x_n^{i_n}\in\TT^n$, we will use the symbol
$D(t)$ as $D(i_1,\dots,i_n)$. Moreover, we write ${\cal
D}=\{D(t)|t\in\TT^n\}$ and denote by Span$_\ccc({\cal D})$ the
$\ccc$-vector space generated by ${\cal D}$. We now introduce some morphisms on ${\cal
D}$ that act as "derivative" and "integral":
\begin{equation}\label{sigma}
\sigma_{x_j}\,(D(i_1,\dots,i_n))\,=\,
\left\{
\begin{array}{ll} D(i_1,\dots,i_j-1,\dots,i_n) & {\rm if}\,i_j>0 \cr 0 & {\rm
otherwise} \cr
\end{array}
\right.
\end{equation}
\begin{equation}\label{ro}
\rho_{x_j}\,(D(i_1,\dots,i_n))\,=\,D(i_1,\dots,i_j+1,\dots,i_n)
\end{equation}
Such operators extend trivially on Span$_\ccc({\cal D})$ by
linearity, and one can easily define $\sigma_t$ and $\rho_t$ for any
$t\in\TT^n$ by composition.
\begin{definition}
A subspace $L$ of Span$_\ccc({\cal D})$ is said to be {\it closed} if
$$\sigma_{x_j}(L)\subseteq L,\,{\it for}\,{\it all}\, j=1,\dots,n.$$
\end{definition}
\begin{definition}
Let $I$ be a primary ideal in $R$ such that
$\sqrt{I}=(x_1,\dots,x_n)$. We define the subspace of differential
operators associated to $I$ as
$$\Delta(I):=\{L\in{\rm Span}_\ccc({\cal D})\,|\,L(f)(0,\dots,0)=0\,{\it for}\,{\it all}\, f\in I\}.$$
Similarly, we associate to each subset $V\subseteq{\rm
Span}_\ccc({\cal D})$ an ideal
$${\mathcal I}(V):=\{f\in R\,|\,L(f)(0,\dots,0)=0\,{\it for}\,{\it all}\, L\in
V\}$$
\end{definition}
\begin{theorem}\label{mmmcorrespondence}
Let ${\frak m}$ be the maximal ideal $(x_1,\dots,x_n)$ of $R$.
There is a bijective correspondence between ${\frak m}$-primary
ideals of $R$ and closed subspaces of Span$_\ccc({\cal D})$
$$\{{\it {\mathfrak m}-primary\; ideals\; in\;} R
\}\TToupdown{\Delta}{\cal I}\{{\it closed\; subspaces\; of}{\rm\;
Span}_\ccc({\cal D}) \}$$ so that $I={\cal I}\Delta(V)$ and
$V=\Delta{\cal I}(I)$ for every $I$ and $V$. Moreover, for a
zerodimensional ${\frak m}$-primary ideal of $R$ whose
multiplicity is $\mu$, we have that ${\rm
dim}_\ccc(\Delta(I))=\mu$.
\end{theorem}
Theorem \ref{mmmcorrespondence} shows that the noetherian operators
associated to a zerodimensional primary ideal form a closed subspace
of Span$_\ccc({\cal D})$. In addition, when considering a
zerodimensional primary ideal, since the dimension of $\Delta(I)$ is
finite, we can view a basis of $\Delta(I)$ as a set of noetherian
operators which, {in this particular case}, happen to be operators
with constant coefficients. Moreover, such a vector space has the
nice property of being closed, fact that has been used by the
authors of \cite{mmm} to construct a
procedure that, given $I$, computes $\Delta(I)$. The algorithm is described below. \\
%\newpage
\begin{algorithm}\label{mmmalgorithm}
Let $I$ be a zerodimensional
primary ideal of $R$ such that $V(I)=\{(0,\dots,0)\}$and let $\mu={\rm dim}_{\ccc}(R/I)$ be its multiplicity.
The following procedure computes the noetherian operators associated to $I$:\\
\begin{quote}
\par\noindent{\bf Input}: ${\cal{G}}=\{g_1,\ldots ,g_t\}$ a \GB for I.
\par\noindent {\bf Output}: $\Delta(I)=\{L_0,\ldots ,L_{\mu-1}\}$
\par\noindent {\bf Initialization}: $i=1$, $L_0=1=$Id$_{{\rm Span}_\ccc({\cal D})}$
\par\noindent {\bf If} $\mu>1$, construct a linear operator $L_1=\sum_{j=1}^n c_j\pp {x_j}$ with an opportune choice
\par\noindent  of the $c_j$'s such that $L_1(f)(0,\dots,0)=0$ is satisfied for each generator $f$ of $I$.
\par\noindent Put $i=2$.
\par\noindent {\bf While} $i<\mu$ {\bf do}
\begin{quote}
\par\noindent define $L_{i+1}$ as a linear combination of $\rho_{j_0}(L_0),\dots,\rho_{j_i}(L_i)$ such that
\par\noindent - $\langle L_0,\dots,L_{i+1}\rangle$ is closed and
\par\noindent - $L_{i+1}(f)(0)=0$ for each generator $f$ of $I$
\end{quote}
\end{quote}
\end{algorithm}
\noindent
\begin{corollary}\label{mmmcorolloperdegrees}
Let $L$ be an operator of $\Delta(I)$, where $I$ is as in
algorithm \ref{mmmalgorithm} and $\mu$ is its multiplicity. Then
{\rm deg}$(L)<\mu$ as an element of $A_n$.
\end{corollary}
\begin{proof}
The construction of $\Delta(I)$ starts with $L_0=1$ and at each
step the degree of $L_{i+1}$ increases of at most $1$, so that the
last element $L_{\mu-1}$ has degree at most $\mu$.
\end{proof}

\begin{remark}{\rm
Algorithm \ref{mmmalgorithm} consists basically in the solution of a
system of linear equations in the coefficients $c_j$ of the linear
combinations $L_{i+1}=c_0 \rho_{j_0} L_0+\cdots +c_i\rho_{j_i} L_i$.
Since the system can have more than one solution, one may simply
pick the one with minimal norm. An implementation for a simplified
version of \ref{mmmalgorithm} has been coded
for \cocoa and is available through the \coala webpage \cite{coala}.\\
%\begin{verbatim}
%http://www.tlc185.com/coala
%\end{verbatim}
}\end{remark}

\begin{example}\label{mmmparabolaexample} {\rm
The following example is taken from \cite{ehrenpreis} (p. 37, ex.
4). Here we show how to study it using algorithm \ref{mmmalgorithm}.
Let us consider the primary ideal at the origin
$I=(y^2,x^2-y)\subset \ccc[x,y]$ whose multiplicity is $4$. We start
with $L_0=1$ and an obvious choice for a linear operator is $L_1=\pp
{x}$. This has also a geometric interpretation: the origin is the
intersection of the two curves given by the generators $y^2$ (the
$x$-axis twice) and $x^2-y$ (a parabola with vertex at the origin).
Such two curves not only intersect at the origin but they are also
tangent along the direction of the $x$--axis, therefore $L_1=\pp
{x}$ must be a noetherian operator. The higher degree operators
describe a higher contact of the line and the parabola at zero. We
can try to find the next one as a combination $L_2=a\pp {x}+b\pp
{xy}$. However, this operator $L_2$ does not respect the closure
condition since $\sigma_x(L_2)=a+b\pp {y}$ which is not in the
subspace $\langle L_0,L_1\rangle =\langle 1,\pp {x}\rangle$. A
different choice for the morphisms $\rho_{x_j}$, instead, gives
$L_2=a\rho_y(1)+b\rho_x(\pp x)=a\pp {y}+b\pp {x^2}$ which respects
closure and annihilates the generators of $I$ at zero with $a=1$ and
$b=\frac{1}{2}$. Again, this operator could have been foreseen in
advance since it is the global annihilator of $x^2-y$ and it
annihilates $y^2$ at the origin. As a last operator, one can choose
$L_3=\rho_x(L_2)=\pp {xy}+\frac{1}{6}\pp {x^3}$. Of course, the
choice $\rho_y(L_2)=\frac{1}{2}\pp {y^2}+\frac{1}{2}\pp {x^2y}$
would have been possible as far as the annihilation of $I$ is
concerned, but it would have violated closure since
$\sigma_x(L_2)=\pp {xy}$ is not a combination of the previous
operators. The iteration ends here since we have found $4$
differential operators. }\end{example}

\subsection{Forward reduction}

We are now going to present an alternative procedure to compute the
noetherian operators associated to $I$ that makes no use of linear
algebra and utilizes the power of \eGGBB.

\begin{algo}[Computation of noetherian operators for zerodimensional ideals]\label{algotaylor}

Let $I$ be a zerodimensional
primary ideal of $R$ such that $\V(I)=\{(0,\dots,0)\}$. The following procedure computes the noetherian operators associated to $I$:\\
\begin{quote}
{\bf Input}: ${\mathcal G}=\{g_1,\ldots ,g_t\}$ a \GB for I.
\par\noindent {\bf Output}: $\Delta(I)=\{L_1,\ldots ,L_\mu\}$.
\begin{quote}
\par\noindent $\bullet$ Compute $\mu(I)={\rm dim}_{\ccc}(R/I)$.
\par\noindent $\bullet$ Write the Taylor expansion at the origin of a polynomial $h\in R$
\par\noindent up to the degree $\mu -1$ with coefficients $c_{\alpha}\in\ccc$:
$$T_{\mu-1}h(x_1,\dots,x_n)=\sum_{\alpha\in\nn^n,\,|\alpha|<\mu}c_{\alpha}x_1^{\alpha_1}\dots
x_n^{\alpha_n}$$
\par\noindent $\bullet$ Write the Normal Form of $T_{\mu-1}h$ with respect to $\cal{G}$ as
\begin{equation}\label{nfideal}
 {\NF}T_{\mu-1}h(x_1,\dots,x_n)=
\sum_{\beta}d_{\beta}x_1^{\beta_1}\dots x_n^{\beta_n}
\end{equation}
\par\noindent and find scalars $a_{\beta\alpha}\in\ccc$ such that $d_\beta=\sum_{\alpha}a_{\beta\alpha}c_\alpha$.
\par\noindent $\bullet$ For each $\beta$ such that $d_\beta\neq 0$, return the operator
$$
 L_\beta=\sum_{\alpha}a_{\beta\alpha}\frac{1}{\alpha_1 !\cdots \alpha_n !} {\pp x_1^{\alpha_1}\cdots\pp
 x_n^{\alpha_n}}=\sum_{\alpha}a_{\beta\alpha}D({\alpha_1}\ldots{\alpha_n}).
$$
\end{quote}
\end{quote}
%%
%
%
%
%DI NUOVO DEVO METTERLO IN FORMA BELLA Let $I$ be a zerodimensional
%primary ideal of $R$ such that $V(I)=\{(0,\dots,0)\}$. Moreover,
%let $\sigma$ be any term ordering on $\TT^n$. Consider the
%following list of
%instructions:\\
%{\rm
%1) Calculate the multiplicity of $R/I$, $\mu(I)={\rm
%im}_{\ccc}(R/I)$\\
%2) Write the Taylor expansion at the origin of a polynomial $h\in
%R$ up to the degree $\mu -1$ with coefficients
%$c_{\alpha}\in\ccc$:
%$$T_{\mu-1}h(x_1,\dots,x_n)=\sum_{|\alpha|=0,\,\alpha\in\nn^n}^{|\alpha|=\mu-1}c_{\alpha}x_1^{\alpha_1}\dots
%x_n^{\alpha_n}$$\\
%3) Calculate the $\sigma$-\GB $\cal{G}$ of $I$. Questo diventa l'input\\
%4) Write the Normal Form of $T_{\mu-1}h$ with respect to $\cal{G}$
%as
%\begin{equation}\label{nf}
% {\NF}T_{\mu-1}h(x_1,\dots,x_n)=
%\sum_{\beta}L_{\beta}x_1^{\beta_1}\dots x_n^{\beta_n}
%\end{equation}
%and find scalars $a_{\beta\alpha}\in\ccc$ such that
%$$d_\beta=\sum_{\alpha}a_{\beta\alpha}c_\alpha$$
%5) For each $\beta$ such that $d_\beta\neq 0$, return the operator
%\begin{equation}\label{zerodimops}
% L_\beta=\sum_{\alpha}a_{\beta\alpha}\frac{1}{\alpha_1 !\cdots \alpha_n !} {\pp x_1^{\alpha_1}\cdots\pp
% x_n^{\alpha_n}}=\sum_{\alpha}a_{\beta\alpha}D({\alpha_1}\ldots{\alpha_n}).
%\end{equation}
%}
%\noindent The above list is an algorithm that provides the
%noetherian operators relative to the ideal $I$.
\end{algo}
\begin{proof}
Let $h(x_1,\dots,x_n)=\sum_{|\alpha|=0}^{
deg(h)}c_{\alpha}x_1^{\alpha_1}\dots x_n^{\alpha_n}$ be the Taylor
expansion centered at the origin of a polynomial $h\in R$ and let
${\cal G}$ be the \GB of $I$. From the theory of \GGBB we know that
the normal form with respect to ${\cal G}$ of $h$ is zero if and
only if $h\in I$, so the condition ${\NF}h=0$ is the one that we
want to characterize. It suffices to write
\begin{equation}\label{finalnf}
\NF (\sum_{|\alpha|=0}^{deg(h)}c_{\alpha}x_1^{\alpha_1}\dots
x_n^{\alpha_n})=\sum_{|\beta|=0}^{deg(h)}d_{\beta}x_1^{\beta_1}\dots
x_n^{\beta_n}=0
\end{equation}
 and deduce from the annihilation of each
coefficient $d_\beta$ in (\ref{finalnf}) a differential condition on
$h$. This completely characterizes the membership of a polynomial
$h$ to $I$. The only thing to observe is that we do not need to work
with terms up to $deg(h)$ for the Taylor expansion. In fact, the
number of differential conditions we need is precisely $\mu$, and so
from corollary \ref{mmmcorolloperdegrees} it follows that the
derivatives to be considered are, in the worst case, the ones of
order $\mu-1$ (see also \cite{mourrain}). Those differential
conditions arise by using coefficients $c_\alpha$ up to
$|\alpha|=\mu-1$. Therefore the Taylor expansion can be truncated at
$\mu-1$.
\end{proof}
\begin{remark}{\rm
It is crucial to observe that we do not need to characterize the
membership of a polynomial $h$ of undetermined degree $deg(h)$ since
we have the bound $\mu-1$ on its degree. Thus algorithm
\ref{algotaylor} is a procedure that is implementable on any
computer algebra software package. Moreover, the computation of the
normal form (\ref{nfideal}) can be done degree by degree, so that we
can stop the reduction process whenever the normal form of a
particular degree is zero. This actually speeds up the computations
in most cases we studied (up to date CPU times for several example
are available on \cite{coala}). }\end{remark}

\begin{example}\label{taylorparabolaexample}{\rm
Consider again Example \ref{mmmparabolaexample} to show the
substantial difference between procedures \ref{mmmalgorithm} and
\ref{algotaylor}. Since $\mu(I)=4$ we start by writing the
truncated Taylor expansion of a polynomial $h\in\ccc[x,y]$:
$$T_3h(x,y)=c_{00}+c_{10}x+c_{01}y+c_{20}x^2+c_{11}xy+c_{02}y^2+c_{30}x^3+c_{21}x^2y+c_{12}xy^2+c_{03}y^3$$
and perform the normal form computation using $x^2\rightarrow y$ and
$y^2\rightarrow 0$ as rewrite rules. Grouping like terms we can
write the remainder of $T_3h$ as a linear combination of the
generators $1,x,y,xy$ of $R/I$ as follows:
\begin{equation}\label{taylorparabolaexample_equation1}
\NF (T_3h)=[c_{00}]+[c_{10}]x+[c_{01}+c_{20}]y+[c_{11}+c_{30}]xy
\end{equation}
We call these four terms a {\it Macaulay basis} for the ideal I,
although this name is also used by some authors for a generalization
of a \eGB. Note that the terms $y^2,x^2y,xy^2$ and $y^3$ disappeared
since they all rewrote to zero. The computation ends by expressing
the coefficients written into square brackets in
(\ref{taylorparabolaexample_equation1}) as operators according to
their meaning as Taylor coefficients. Namely $[c_{00}]\rightarrow 1,
[c_{10}]\rightarrow \pp {x}, [c_{01}+c_{20}]\rightarrow \pp
{y}+\frac{1}{2}\pp {x^2}, [c_{11}+c_{30}]\rightarrow \pp
{xy}+\frac{1}{6}\pp {x^3}$. This gives the same result obtained in
the example \ref{mmmparabolaexample} as expected. This is not
surprising since theorem \ref{mmmcorrespondence} states that the
correspondence $I\leftrightarrow\Delta(I)$ is one-to-one.
}\end{example}

Algorithm \ref{algotaylor} does not take directly into account the
closure of the space of noetherian operators, as algorithm
\ref{mmmalgorithm} did. The fact that $\Delta(I)$ is closed is a
general fact which follows from a Leibniz formula for the morphisms
$\sigma_{x_j}$ and the fact that $I$ is an ideal (see \cite{mmm2},
prop. 2.4). This is true not only for zerodimensional ideals but
also in positive dimension, as we will see in section $9$. We want
to show that the closure of $\Delta(I)$ is also a direct consequence
of algorithm \ref{algotaylor} and of the following property of
Macaulay bases.
\begin{lemma}
Let $I\subset R$ be an ideal and let ${\mathcal M}$ be the
Macaulay basis of $R/I$, i.e. the generators of $R/I$ as a
$\ccc$-vector space. Let $s_{x_j}:\TT^n\rightarrow\TT^n$ be the
"derivative" morphism
\begin{equation}\label{s_morphism}
s_{j}(x_1^{i_1}\cdots x_n^{i_n})= \left\{
\begin{array}{ll}
 x_1^{i_1}\cdots x_j^{i_j-1}\cdots x^n_{i_n} & {\it if\def\terms#1{{\rm \TT^{#1}}}}\ i_j>0 \cr
 0 & {\it otherwise} \cr
\end{array}
\right.
\end{equation}
Then ${\mathcal M}$ is $s_j$-closed for each $j$.
\end{lemma}
\begin{proof}
It is known that the Macaulay basis for $R/I$ can be computed
through a \GB ${\cal G}$ of $I$. In fact it is (see \cite{kr},
theorem 1.5.7):
$${\mathcal M}=\TT^n\backslash \LT_\sigma({\cal G})$$
where $\sigma$ is any term ordering on $\TT^n$. Since ${\cal G}$
is a \GB for $I$, the leading term ideal $\LT_\sigma(I)$ coincides
with $\LT_\sigma({\cal G})$. Let $t\neq 0$ be a term of ${\mathcal
M}$. Suppose that there exists an index $j$ such that $0\neq
s_j(t)\notin {\mathcal M}$. Then $s_j(t)\in \LT_\sigma({\cal G})$.
The latter being an ideal, we have $t=x_j\cdot s_j(t)\in
\LT_\sigma({\cal G})$, which is a contradiction. Note that if $
s_j(t) \notin {\mathcal M}$ for all $j$, this simply says that $t=0$
which is again a contradiction.
\end{proof}
The morphism $s_j$ introduced in the above lemma is the analogue
of $\sigma_{x_j}$ defined in section \ref{closure}, and we
will show in the next proposition that the $s_j$-closure of
${\mathcal M}$ is equivalent to the $\sigma_{x_j}$-closure of the
space of noetherian operators associated to $I$.

\begin{prop}\label{closureprop}
Let $I$ be a zerodimensional primary ideal of $R$ such that
$\V(I)=\{(0,\dots,0)\}$ and let ${\cal O}=\{L_\beta\}$ be the set of
operators computed with algorithm \ref{algotaylor}. Then
Span$_\ccc(\{L_\beta\})$ is a closed subspace of Span$_\ccc({\cal
D})$.
\end{prop}
\begin{proof}
Let $L_\beta\in{\cal O}$, and let $d_\beta$ be the corresponding
coefficient of the normal form $\NF(h)$ as computed with the algorithm. Let $x^\beta=x_1^{\beta_1}\cdots x_n^{\beta_n}$ be
the term whose coefficient is $d_\beta$. It is clear that such a
term is part of the Macaulay basis of $R/I$ since it appears in
the expression of $\NF(h)$, which is a representation of the class
of $h$ in the quotient $R/I$. Denote by $F_{\beta}$ the set of
operators of ${\cal O}$ such that the corresponding term in the
expression of $\NF(h)$ divides $x^\beta$:
$$F_\beta=\{L_\gamma\in{\cal O}\,{\rm such}\,{\rm
that}\,x^\gamma|x^\beta\}$$ and for each $L_\gamma\in F_\beta$
consider $t_\gamma=x^{\beta-\gamma}$. Since each $L_\gamma$ has been
computed from the Taylor expansion of using a division algorithm
that uses a \GB ${\cal G}$ of $I$, we have that (see \cite{kr},
prop. 2.2.2) if $h'$ is such that
$$x^\gamma=\NF(h')\;{\rm and}\;{\rm supp}(h')\subseteq{\rm supp}(h)$$
 then $$x^\beta=t_\gamma x^\gamma=\NF(t_\gamma)
\NF(h')=\NF(t_\gamma h')$$
 i.e. the term in $x^\beta$ is obtained rewriting a multiple of
 that part of the polynomial $h$ which rewrites to $x^\gamma$. By looking at
 the expression of $L_\beta$ is then  obvious that
 $$ \sigma_{t_\gamma} (L_\beta)=L_\gamma$$
 since $L_\beta$ is written as a combination of Taylor coefficients
 corresponding to the terms of $t_\gamma h'$. It now suffices to
 prove that such $t_\gamma$'s are enough to conclude that ${\cal
 O}$ is closed. This is a consequence of the
 previous lemma, since all the $d_\gamma$ in $F_\beta$ are
 associated to those terms $x^\gamma$ of the Macaulay basis ${\mathcal M}$
 that divide $x^\beta$, hence from the $s_j$-closure of ${\mathcal M}$ we deduce that
$\{x^\gamma\,=s_{t_\gamma}(x^{\beta})\}=\{s_j(x^{\beta}),j=1,\dots,n\}$.
\end{proof}
\subsection{Backward reduction}
We could think of performing the reduction step of the algorithm for
the computation of noetherian operators for zerodimensional ideals
"backwards". Instead of writing the full Taylor expansion and then
using the \GB of $I$ to rewrite it, we start from the residual
monomials, which are easily calculated for example with \ecocoa. We
then "pull back" each monomial using the generators of $I$ as
"anti-rewrite rules". Let us explain what we mean by this. In
general, when using a polynomial $f$ to rewrite another polynomial
$g$, we use its leading monomial $\LT(f)$ to divide the polynomials
$g$ and then we substitute each $\LT(f)$ in $g$ with the {\it tail}
of $f$, $\LT(f)-f$. For instance, we rewrite $g=x^3$ to $xy$ using
$f=x^2-y$, by replacing $x^2$ in $x^3$ with the tail
$x^2-(x^2-y)=y$. This operation, when performed using the elements
of a \GB for $I$, does not alter the class of $g$ in $R/I$ and leads
to the normal form $\NF(g)$. What we mean by "anti-rewriting" is,
roughly speaking, to use the smallest monomial of $f$, $in(f)$, and
replace it with the head of the polynomial, $in(f)-f$. This way,
from $in(f)$ we "climb up" to find all the other monomials that are
equivalent to $in(f)$ modulo $(f)$. Here is a more precise
definition.
\begin{definition}\label{backwards-rewrite}
Let $f$ be a polynomials of $R$, let $g$ be a monomial and let
$m=in(f)$ be the smallest term of $f$ with respect to a given term
ordering on $\TT^n$. We say that $g$ {\it rewrites backwards} to
$g'$ in one step, using $f$, if $m$ divides $g$ and
$$g'=\frac{g}{m}(m-f).$$
\end{definition}
\begin{example}{\rm
With this terminology, $g=xy$ rewrites backwards to $x^3$ using
$x^2-y$, which is exactly the opposite of the standard rewrite
process that leads from $x^3$ to $xy$. If we use $f=x^2+xy-2y$
instead, $g=xy$ rewrites to $g'=\frac{1}{2}x^3+\frac{1}{2}x^2y$.
Finally, $g$ could not be rewritten backwards suing $x^2-y^2$ since
$y^2$ does not divide $g$. Notice that in general if we perform a
one-step backward reduction and then a one-step reduction in the
usual way, we obtain back $g$.}
\end{example}
We can now apply an iteration of this procedure of rewriting
backwards a monomial using a \GB for $I$. We start from a residual
monomial and we rewrite it backwards using one generator. Then we
rewrite backwards each monomial obtained after this step, if
possible, using any element of the \eGB. Technically this procedure
never ends, as we can imagine, to obtain a new polynomial of higher
degree at each step, as for example with $g=x$ and $f=x^2-x$.
However, for the purpose of computing noetherian operators, we know
from section \ref{zerodimsection} that, as polynomials in $\ccc[\pp
x_1,\dots,\pp x_n]$, they have degree at most $\mu-1$. Therefore we
can stop the iteration once we have reached a polynomial of such
a degree. Let us illustrate this idea with an example before we
present the algorithm in general.
\begin{example}\label{backwards_noetherop_example}{\rm
Consider the ideal $J=(x^2-z,y^2-z,z^2)$ in $\ccc[x,y,z]$. It
represents the origin in $\ccc^3$ with multiplicity eight. Its
generators are a {\tt DegLex} \eGB. The residual monomials for $R/J$
are
$$\{1,x,y,xy,z,xz,yz,xyz\}.$$
First, let us reconstruct the noetherian operators associated to
$xyz$. By rewriting it using $x^2-z$ we obtain the new monomial
$x^3y$. This cannot be rewritten further. However, the term $xy^3$
is another monomial that is "attracted" by $xyz$ via the other
generator $y^2-z$ of $J$. Summing up the residual monomial and all
the results of the backward reduction we then obtain
$g'=x^3y+xy^3+xyz$
 whose dual $D(g')=\frac{1}{6}\pp x^3\pp y+\frac{1}{6}\pp x\pp y^3+\pp x\pp y\pp
 z$ is actually  the noetherian operator of $J$ relative to $xyz$.}
\end{example}
The choice of the residual monomial $xyz$ in Example
\ref{backwards_noetherop_example} is not random. Indeed it is
maximal among all the residual monomials with respect to the
derivative morphisms (\ref{s_morphism}).
\begin{definition}\label{cornermonomial}
Let $m$ be a residual monomial of $R/I$. We say that $m$ is a {\it
corner monomial} if it is maximal with respect to the monoid
structure of $\TT^n$, i.e. if
$$x_i\cdot m\in\LT(I),\;{\it for\;all\;}i=1\dots n.$$
\end{definition}
If we represent $R/I$ as a subset of $\nn^n$, the corner monomials
are exactly in corner position. Proposition \ref{closureprop} says
that the noetherian operators are generated by the ones
corresponding to the corner monomials by taking the closure with
respect to the morphisms (\ref{sigma}). This fact allows to come up
with a general procedure that constructs the noetherian operators
starting with the corner monomials and then generates the entire
space of noetherian operators.
\begin{algorithm}\label{backwards_algo}{\it
Let $I\subset R$ be a zerodimensional primary ideal of multiplicity
$\mu$ centered at the origin. The following list of instructions
construct the noetherian operators associated to $I$:\\
\begin{quote}
\par\noindent{\bf Input}: a \GB $\G$ of $I$ and the residual monomials of $R/I$.
\par\noindent {\bf Output}: the space of noetherian operators associated to $I$.
\par\noindent $\bullet$ Construct the set $\cc$ of corner monomials using definition \ref{cornermonomial}.
\par\noindent $\bullet$ For each corner monomial $m\in\cc$ find the associated noetherian operators by rewriting
 it backwards with respect to $\G$ using
definition \ref{backwards-rewrite}. Stop when the backward reduction
is not possible anymore or when the degree of the polynomial
obtained is $\mu-1$.
\par\noindent $\bullet$ Collect all the polynomials obtained in the set ${\cal D}$.
\par\noindent $\bullet$ Compute the closure of ${\cal D}$ by applying the morphism $\sigma_{x_i}$, $i=1\dots n$ to all its elements.
\par\noindent $\bullet$ For each element $L$ in the closure of ${\cal D}$
calculate $D(L)$.
\end{quote}}
\end{algorithm}
\subsection{Extension to modules}
All the results in the previous subsections can be extended in a straightforward fashion to the case of zerodimensional primary modules. Rather than giving the details, we use the \cocoa version of the algorithm for modules to look at a couple of examples.

\begin{example}{\rm
Let $A$ be the matrix
$$A=
\left(
\begin{array}{cc}
x & 1 \cr y & x \cr 0 & y \cr
\end{array}
\right)
$$
and let $M$ be the module generated by the rows of $A$, i.e.
$M=\langle xe_1+e_2,\,xe_2+ye_1,\,ye_2\rangle$. The module term
ordering we choose is {\tt Lex-Pos}, meaning that to compare two
terms we first look at the power product, using ${\tt Lex}$, and
then we look at the position. The way we just wrote the generators
of $M$ reflects this choice. It is clear that $J_M=(x^2-y,y^2,xy)$,
and, using for example \ecocoa, we find out
that:\\
- $\mu(M)=3$ \\
- the ${\tt Lex}$--\GB of $M$ is ${\cal
G}=\{xe_1+e_2,\,xe_2+ye_1,\,ye_2,\,y^2e_1\}$\\
- a Macaulay basis for $M$ is the set $\{e_1,e_2,ye_1\}$.\\
We begin by writing explicitly the vectorial Taylor expansion of a
vector $w(x,y)\in R^s$ up to degree 2:\\
$T_2w(x,y)=
 { c^1_{00}e_1+c^2_{00}e_2+c^1_{10}xe_1+c^2_{10}xe_2+c^1_{01}ye_1}
+c^2_{01}ye_2+ {c^1_{20}x^2e_1}
+c^2_{20}x^2e_2+c^1_{11}xye_1+c^2_{11}xye_2+c^1_{02}y^2e_1+c^2_{02}y^2e_2$.\\
Only few terms survive after we compute the normal form relative to
the \GB ${\cal G}$, leading to
$$\NF(w)=[c^1_{00}]e_1+[c^2_{00}-c^1_{10}]e_2+[c^1_{20}+c^1_{01}-c^2_{10}]ye_1.$$
We conclude that the noetherian operators associated to $M$,
written in vectorial form, are
$$D^1_{00}=(1,0),\quad D^2_{00}=(-\pp x,1),\quad
D^1_{01}=(\frac{1}{2}\pp {x^2}+\pp y,-\pp x)$$ and it is easy to
check that they generate a closed subspace since
$\sigma_x(D^2_{00})=\sigma_y(D^1_{01})=D^1_{00}$ and
$\sigma_x(D^1_{01})=D^2_{00}$.
 }
\end{example}

\begin{example}[Solution of a system of PDEs]\label{pdesystem}
{\rm In the introduction we saw that the Fundamental Principle can be
used to write an integral representation of the solution of a
system of linear constant coefficient partial differential
equations. We will show how this can be applied, now that we know
how to compute noetherian operators. Consider the overdetermined PDE system given by
\begin{equation}\label{thesystem}
\left\{
\begin{array}{rcr}
f_{zz}-f_z+f_t+2g_z  &   =   &  g  \cr
f_{zt}+g_t           &   =   &  0  \cr
f_{tt}+g_{zt}-g_t    &   =   &  0  \cr
f_{t}-g_{zz}+g_z+g_t &   =   &  0  \cr
\end{array}
\right.
\end{equation}
where $f,g\in C^{\infty}(\rr^2)$ and we use indices to denote derivatives.
The general solution to (\ref{thesystem}) can be written using a generalization of (\ref{3.3}). We consider the
rectangular operator $P(D)$ defined by
$$ P=
\left(
\begin{array}{cc}
x^2-x+y & 2x-1 \cr xy & y \cr y^2 & xy-y \cr y & x^2-x-y \cr
\end{array}
\right)
$$
where $x$ and $y$ are the dual variables of $z$ and $t$
respectively. Note that we are choosing a particular Fourier
transform to write $P(D)$ so that it does not take into account the
factor $-\sqrt{-1}$. The module $M$ associated to the matrix $P$ is
not primary, hence we can use Singular to get a primary
decomposition (using the function {\tt modDec} form the library {\tt
mprimdec.lib}). $M$ is the intersection of the two
zerodimensional modules\\
$M_1=\langle (x,1),(y,x),(0,y)\rangle,\quad J_1=\sqrt{M_1}=(x,y)$\\
$M_2=\langle (x-1,1),(y,0),(y,x-1)\rangle,\quad J_2=\sqrt{M_2}=(x-1,y)$\\
of multiplicity, respectively, $3$ and $2$. We already computed the
operators associated to the module $M_1$ in the previous example. To
compute the operators associated to $M_2$ we need to shift the
variety to the origin using the change of coordinates $(X=x-1,
Y=y)$. Then, using the new variables $X$ and $Y$, we can apply the
module version of algorithm \ref{algotaylor} and find the noetherian
operators: $\{(1,0),(\pp X,-1)\}$. Going back to the variables $x,y$
we have the set $\{(1,0),(\pp x,-1)\}$. Therefore, it is possible to
write explicitly the solutions to (\ref{thesystem}) as

%\begin{equation}
$$
\left(
\begin{array}{c}
f(z,t) \cr g(z,t) \cr
\end{array}
\right)=
A\left(
\begin{array}{c}
1 \cr 0 \cr
\end{array}
\right)_{|_{(0,0)}} e^{zx+ty}
+B \left(
\begin{array}{c}
-\pp x \cr 1 \cr
\end{array}
\right)_{|_{(0,0)}} e^{zx+ty}
+C \left(
\begin{array}{c}
\frac{1}{2}\pp {x^2}+\pp y \cr -\pp x \cr
\end{array}
\right)_{|_{(0,0)}} e^{zx+ty} +
%\end{equation}
$$
$$
%\begin{equation}
+
D
\left(
\begin{array}{c}
1 \cr 0 \cr
\end{array}
\right)_{|_{(1,0)}} e^{zx+ty}
+E \left(
\begin{array}{c}
\pp x \cr -1 \cr
\end{array}
\right)_{|_{(1,0)}} e^{zx+ty} =
\left(
\begin{array}{c}
 A-Bz+\frac{1}{2}Cz^2+Ct+De^z+Eze^z
\cr B-Cz-Ee^z \cr
\end{array}
\right)
%\end{equation}
$$
}
\end{example}

\section{The case of positive dimension}

When dealing with ideals and modules whose dimension is positive, in general
one may not expect the associated noetherian operators to be
constant coefficient linear operators. In fact, this is
the case for some of the examples from the literature (see
\cite{ehrenpreis,palamodov}). For instance, when considering the
ideal $I=(x^2,y^2,-xz+y)\subset\ccc[x,y,z]$ one has that a set of
noetherian operators associated to $I$ is $\{1,\pp x+z\pp y\}$ and
it can be proved that there exist no set of noetherian operators
with constant coefficients associated to $I$ (see \cite{palamodov},
example 4, p. 183). However, an interesting property that we notice
in this case is that the set of "differential" variables from the
set of variables appearing in the polynomial coefficients (in this
case such sets are respectively $\{x,y\}$ and $\{z\}$). This is
actually valid whenever we can put the algebraic variety in a
particular position, through an opportune change of coordinates,
called {\it normal position}. To do this, one can apply the
procedure of Noether normalization to the ideal $I$. This algorithm
comes from the so-called Noether Normalization Theorem (see
\cite{bjo}, p. 116). We now state a version of the theorem that we
will need for our computations:

\begin{theorem}[Noether Normalization Theorem]\label{normalization}
Let ${I}$ be a primary ideal of $\ccc[z_1,\dots,z_n]$. There exist a non--negative integer
$d$ and a (linear) change of coordinates
$$\varphi:\ccc[z_1,\dots,z_n]\rightarrow\ccc[x_1,\dots,x_{n-d},t_1,\dots,t_d]$$ such that:\\
{\rm a)} $\varphi({I})\cap\ccc[t_1,\dots,t_d]=(0)$,\\
{\rm b)} $\ccc[z_1,\dots,z_n]/{I}$ is a finitely generated $\ccc[t_1,\dots,t_d]$--module,\\
{\rm c)} for each $i=1\dots n-d$, $\varphi({I})$ contains a polynomial of the form
$$Q_i(t_1,\dots,t_d,x_i)=x_i^{e_i}+p_1(t_1,\dots,t_d)x_i^{e_i-1}+\cdots+p_{e_i}(t_1,\dots,t_d)$$
 where $e_i$ is the degree of the polynomial $Q_i$. \\
The ideal $\varphi({I})$ is said to be in {\rm normal position} with respect to the
variables $x_1,\dots,x_{n-d}$.
\end{theorem}
\begin{remark}{\rm
The proof of the Normalization Theorem can be found for example in
\cite{bjo}, in the case of prime ideals. However, as shown in
\cite{singularbook}, the result holds for the general case with the
exception of condition a) which requires $I$ to be primary. If the
ideal $I$ is prime, the polynomials $Q_i$ in condition c) can be
chosen to be irreducible. The proof of the theorem provides an
algorithm to achieve the normal position. Basically, at each step
one constructs the polynomial $Q_i$, performing a {\it generic}
coordinate change such that $Q_i$ has a monic leading term of the
form $x_i^{e_i}$, and then one eliminates the variable $x_i$. A
procedure to compute the Noether normalization of an ideal has also
been studied in \cite{logar} and it is available in Singular through
the library {\tt algebra.lib} (see \cite{singular} and its manual).
We coded a version of the algorithm for \cocoa as well,
\cite{coala}. }
\end{remark}
\noindent Theorem \ref{normalization} basically states that it is
possible to find a new system of coordinates where the $x$ variables
act as "variables" and the $t$ variables act as "coordinates", and
where the integer $d$ appearing in \ref{normalization} is nothing
but the dimension of the ideal $I$. Hence, if we make the variables
$t$ invertible, i.e. if we extend the ideal to the ring $\ccc(t)[x]$
where $\ccc(t)$ is the ring of quotients of $\ccc[t]$, we end up
with a zerodimensional ideal. Furthermore, since we are interested
only in primary ideals,  we may expect that the extension of the
ideal to $\ccc(t)[x]$ is still primary. The following proposition
assures that such facts hold if $I$ is in normal position.

\begin{prop}\label{extendedideal}
Let $I=(f_1,\dots,f_r)$ be a primary ideal of dimension $d$ in the
polynomial ring $R=\ccc[x_1,\dots,x_{n-d},t_1,\dots,t_d]$, in
normal position with respect to $x_1,\dots,x_{n-d}$. Denote by
$R_d=\ccc(t_1,\dots,t_d)[x_1,\dots,x_{n-d}]$ the ring of
polynomials in the $x$ variables with coefficients in the field of
fractions
$\ccc(t_1,\dots,t_d)={\rm Frac}(\ccc[t_1,\dots,t_d])$. The following facts hold:\\
{\rm 1)} the inclusion map  $\varphi_{|_I}:\, I\hookrightarrow IR_d$ is
injective and $IR_d\cap R=I$,\\
{\rm 2)} the extended ideal $IR_d$ is primary,\\
{\rm 3)} the extended ideal $IR_d$ is zerodimensional.
\end{prop}
\begin{proof}
The fact that the inclusion is injective is trivial. To prove 1), let us consider a polynomial $f$
in $R\cap IR_d$. As an element of $IR_d$ it can be written in the form
$$f=\sum_{i=1}^{r} \frac{a_i(x,t)}{b_i(t)}f_i(x,t)$$
where $x=(x_1,\dots,x_{n-d})$, $t=(t_1,\dots,t_{d})$, and $a_i$ and
$b_i$ are just polynomials in the set of variables indicated in
parenthesis. Let $b(t)=\prod_{i=1}^rb_i(t)$ and consider the product $bf$.
Both $b$ and $f$ are polynomials in $R$ and their product is an
$R$-linear combination of the generators of $I$, so $bf\in I$. Since
$I$ is primary it follows that either $b^m\in I$ for some positive
integer $m$ or $f\in I$. The first possibility is in contradiction
with condition a) of the Noether normalization, hence $f\in I$. This
proves that $IR_d\cap R\subseteq I$. The opposite inclusion is
trivial, so we conclude that $IR_d\cap R=I$. The same type of
argument can be used to prove that $IR_d$ is primary: consider two
fractions
$$f(x,t)=\frac{a(x,t)}{b(t)},\quad g(x,t)=\frac{c(x,t)}{e(t)}$$ such that $fg\in IR_d$. Then
$(bf)\cdot (eg)$ is a polynomial in $I$ and since $I$ is primary  we
either have $bf\in I$ or $e^mg^m\in I$ for some positive integer
$m$. In the first case, using again that $I$ is primary and using
condition a) of Theorem \ref{normalization}, we get that $f$ is in
$I$. In the second case we have that $g^m$ is in $I$. Therefore
either $f \in IR_d$ or $g^m\in IR_d$. Finally, statement 3) follows
from the theory of the dimension of an ideal, since
$(\bar{t_1},\dots,\bar{t_d})$ is a maximal regular sequence in $R/I$
that reduces to just constants when extending the ideal to $R_d$.
\end{proof}

Before we move on and present an equivalent version of algorithm
\ref{algotaylor} for non zerodimensional ideals, there is still one
more step. Formerly, when treating the zerodimensional case, we
chose to start with a \GB for the ideal $I$, computed with respect
to any term ordering. This is no longer possible if we want to
extend the procedure to the positive dimensional case. In fact,
after we perform the normalization, the variables $t$ play the role
of "constants" once we extend $I$ to $R_d=\ccc(t)[x]$. The following
example illustrates a problem that may occur if we do not choose
carefully the term ordering on $R$.
\begin{example}{\rm
Consider the ideal $I=(x^2-t,xt-1)$ in $\ccc[x,t]$. A ${\tt
DegLex}$--\GB for $I$ (with $x>t$) is given by ${\cal G}=\{{\bf
x^2}-t,{\bf x}{\bf t}-1,{\bf t^2}-x\}$, where the leading term are
highlighted in bold. When we look at such polynomials in $R_d$,
however, we see that the leading terms change, in fact the last
polynomial should better be written as ${\bf -x}+t^2$. Note that in
this case the extended ideal $IR_d$ happens to be the whole ring
$R_d$ since the polynomial $t^3-1$ belongs to $IR_d$, and such
polynomial is a constant in $\ccc(t)[x]$. It is a necessary and
sufficient condition for an ideal to be the whole ring that any \GB
with respect to any ordering contains a constant polynomial, but if
we look at ${\cal G}$ we see that there is no such a constant,
meaning a polynomial only in the variable $t$. Therefore we conclude
that the set ${\cal G}$ does not form a \GB for $IR_d$, with respect
to the ordering ${\tt DegLex}$ restricted to the terms in $x$. If we
choose instead the term ordering ${\tt Lex}$, a \GB for $I$ is given
by ${\cal G}=\{-{\bf x}+t^2,{\bf t^3}-1\}$, and in this case it
contains a polynomial in $t$, making ${\cal G}$ a \GB for $IR_d$ as
well. }
\end{example}

As the example shows, we really want the variables $x$ to be the
main variables with respect to which the \GB needs to be computed.
This can be achieved using ${\tt Lex}$, or any other elimination
ordering with respect to the variables $x$. Lemma \ref{lemmaproduct}
then ensures that after extending the ideal to $R_d$, \GGBB are
preserved. We now have all the ingredients to generalize algorithm
\ref{algotaylor} to the case of an ideal of dimension greater than
zero. As in section \ref{zerodimsection}, we will suppose that a
primary decomposition of the ideal has already been calculated.

\begin{algo}[Noetherian operators for positive dimensional ideals]\label{highertaylor}

Let $d$ be a positive integer, $x=(x_1,\dots,x_{n-d})$ and
$t=(t_1,\dots,t_d)$ be variables and let $\sigma=\sigma_x\cdot
\sigma_t$ be a product ordering. Let $I$ be a primary ideal in
$R=\ccc[x,t]$. Suppose that $I$ is in normal position with respect
to $x$. Moreover, let $IR_d$ be the extended ideal in
$R_d=\ccc(t)[x]$ and suppose
that the characteristic variety of  $IR_d$ in $\ccc(t)^d$ is the origin. The following procedure computes the noetherian operators associated to $I$:\\
\begin{quote}
{\bf Input}: ${\mathcal G}=\{g_1,\ldots ,g_r\}$ a $\sigma$--\GB for I.
\par\noindent {\bf Output}: a set of noetherian operators for $I$.
\begin{quote}
\par\noindent $\bullet$ Compute the multiplicity of the ideal, $\mu(I)$.
\par\noindent $\bullet$ Write the Taylor expansion at the origin of a polynomial $h\in\ccc[x]$
\par\noindent up to the degree $\mu -1$ with variable coefficients $c_{\alpha}$:
\begin{equation}\label{taylorh}
\hat{h}:=T_{\mu-1}h(x_1,\dots,x_{n-d})=\sum_{\alpha\in\nn^{n-d}}^{|\alpha|<\mu}c_{\alpha}x_1^{\alpha_1}\dots
x_{n-d}^{\alpha_{n-d}}
\end{equation}
\par\noindent $\bullet$ Let $x^{a_i}t^{b_i}$ be the leading term of $g_i$ and define $t^{\gamma}:=t^{b_1}\cdots t^{b_r}$
%\par\noindent $\bullet$ Define $\hat{h}:=T_{\mu-1}h$
\par\noindent {\bf Repeat}
\begin{quote}
$\bullet$ Multiply $\hat{h}$ by $t^\gamma$ and compute its normal
form with respect to $\cal{G}$. \par\noindent $\bullet$ Rename that
as $\hat{h}$:
\begin{equation}\label{nfhiger}
 \hat{h}:={\NF}(t^{\gamma}\hat{h})=
\sum_{\beta}d_{\beta}(t)x_1^{\beta_1}\dots x_{n-d}^{\beta_{n-d}}
\end{equation}
\par\noindent
\end{quote}
\par\noindent {\bf Until} the number of nonzero $d_\beta$ is exactly $\mu$.
\par\noindent $\bullet$ For each $\beta$ such that $d_\beta\neq 0$, find polynomials $a_{\beta\alpha}(t)$ such that
$d_\beta(t)=\sum_{\alpha}a_{\beta\alpha}(t)c_\alpha$ and return the
operator
$$
 L_\beta=\sum_{\alpha}a_{\beta\alpha}(t)\frac{1}{\alpha_1 !\cdots \alpha_{n-d} !} {\pp x_1^{\alpha_1}\cdots\pp
 x_{n-d}^{\alpha_{n-d}}}=\sum_{\alpha}a_{\beta\alpha}(t)D({\alpha_1},\ldots,{\alpha_{n-d}},0,\dots,0).
$$
\end{quote}
\end{quote}
\end{algo}
\begin{proof}
Let $h$ be a polynomial of $R$. We want to characterize the
membership of $h$ to $I$. Since we are assuming that $I$ is in
normal position, by condition 1) of Proposition \ref{extendedideal}
this is equivalent to the membership of $h$ to $IR_d$. Since the
latter is a zerodimensional ideal of multiplicity $\mu$, $h\in IR_d$
if and only if the Taylor polynomial  of degree $\mu-1$ of $h$, with
coefficients in $\ccc(t)$, reduces to zero when rewriting it using a
\GB for $IR_d$. This follows from the the same proof as in algorithm
\ref{algotaylor}. By Lemma \ref{lemmaproduct}, a $\sigma_x$--\GB for
$IR_d$ is given by the same elements of the \GB of $I$. This means
that computing a normal form in $I$ and in $IR_d$ is equivalent.
However, when writing the Taylor expansion (\ref{taylorh}), we need
to consider that the coefficients $c_{\alpha}$ also depend on $t$.
In order to be able to perform a one-step reduction, we need each
term in (\ref{taylorh}) to be at least multiplied by $t^\gamma$.
This does not affect the membership of $T_{\mu-1}h$ as a polynomial
in $R_d$ since it is just a multiplication by a constant. Also when
considering the expression (\ref{taylorh}) in $\ccc[x,t]$, the
effect of the multiplication does not change the annihilation of
$\NF(T_{\mu-1}h)$, since obviously
$$\NF(T_{\mu-1}h)=0\Leftrightarrow \NF(t^{\gamma}T_{\mu-1}h)=0.$$
The one-step reduction is then iterated enough times in
(\ref{nfhiger}) until we reach a sufficiently small number of
nonzero terms (namely $\mu$). By what we have proved so far, it is
then clear that at the end of the process the polynomial $\hat{h}$
is exactly the normal form of $T_{\mu-1}h$ as a polynomial in
$R_d$ and hence the annihilation of its coefficients is equivalent
to the condition $h\in IR_d$.
\end{proof}
\begin{remark}{\rm
The main difference with respect to the algorithm for zerodimensional ideals is that, in this case, we do not know if
after just one step of reduction we have achieved the normal form of the polynomial $h(x,t)$, since the
multiplication by $t^\gamma$ could not be enough to assure that $h$ has been rewritten to a sum that runs over
just the Macaulay basis terms for $IR_d$. Multiplying $T_{\mu-1}h$ once by $t^\gamma$ is definitely enough for a one-step reduction of each term of the Taylor expansion. That is, each term is being rewritten using
at most one of the elements of the \eGB. However, further reductions might occur if we multiply again by $t^\gamma$. Also,
note that such an iteration has to terminate because $\sigma_{x}$ is a well ordering.
}
\end{remark}
\begin{remark}{\rm
The reduction step (\ref{nfhiger}) for ideals with few generators is
not very heavy, but performing it multiple times could slow down the
procedure by a significant amount. We believe that it is possible to
find an exponent $\gamma_1$ large enough so that we need to multiply
by $t^{\gamma_1}$ just once, allowing the reduction to bring
$\hat{h}$ all the way down to its final expression. For example,
choosing $\gamma_1=\mu\cdot \gamma$ seems to work fine at least in
the cases we tested, without the need of further iteration.
}\end{remark} When applying algorithm \ref{highertaylor} to an ideal
$I$ in normal position, some redundant factors in $t$ could appear
as an effect of the iterative multiplication by $t^\gamma$ at each
step. Since such factors are constants in $R_d$, they are actually
not needed to characterize the membership of a polynomial in $R_d$.
It is then possible to eliminate these factors from the final
expression of the noetherian operators. The next example will
clarify what we mean.
\begin{example}{\rm
Consider the system of partial differential equations in three variables given by
$$
\left\{
\begin{array}{ccc}
f_{xx} & = & 0 \cr
f_{yy} & = & 0 \cr
f_{y} & = & f_{xt} \cr
\end{array}
\right. .$$ Its solutions are differentiable functions of the form
$f(x,y,t)=A(t)+B(t)x+B'(t)y$, where $A$ and $B$ are arbitrary
functions of $t$. We want to derive this last statement using the
fundamental principle. The primary ideal associated to the system is
$I=(x^2,y^2,-xt+y)$ in $\ccc[x,y,t]$ (see \cite{palamodov}). If we
consider the ${\tt Lex}$ ordering where $x>y>t$, a \GB for $I$ is
 given by $(x^2,xy,y^2,-xt+y)$. Let us compute the associated noetherian operators using algorithm \ref{highertaylor}.
 It is immediate to check that $I$ is in normal position with respect to
$x$ and $y$ and that, after inverting $t$, the variety associated to
$I\ccc(t)[x,y]$  is the origin in $\ccc(t)^2$. The multiplicity of
$I$ can be computed with \ecocoa, and it is $\mu=2$. So we just need
to write a linear polynomial $h$ with variable coefficients and
multiply it by $t$, which is the only term in $t$ appearing in the
leading terms of the \eGB:
$$T_1\hat{h}=t\cdot T_1h=tc_{00}+tc_{10}x+tc_{01}y.$$
The only rewrite rule that we need to use to reduce $h$ is hence $xt\rightarrow y$ which leads to the final
expression for the normal form
$$\NF(\hat{h})=[tc_{00}]+[c_{10}+tc_{01}]y.$$
Since the terms in $x$ and $y$ of the last expression are exactly
$\mu=2$, we do not need to proceed further and then we conclude that
the noetherian operators are $\{t,\pp x+t\pp y\}$. Since the first
is a multiple of $t$, we can divide it by $t$ and get the final set
$\{1,\pp x+t\pp y\}$. Now we can write the integral formula for the
general solution of the system, using $\zeta,\eta,\tau$ as dual
variables:
$$f(x,y,t)=\int\limits_{\zeta=\eta=0}e^{i(x\zeta+y\eta+t\tau)}d\mu_1(\zeta,\eta,\tau)+\int\limits_{\zeta=\eta=0}(\pp \zeta+\tau\pp \eta) e^{i(x\zeta+y\eta+t\tau)}d\mu_2(\zeta,\eta,\tau)=$$
$$=\int_{\rr}e^{it\tau}d\mu_1(\tau)+\int_{\rr}i(x+y\tau)e^{it\tau}d\mu_2(\tau)
=\int_{\rr}e^{it\tau}d\mu_1(\tau)+x\int_{\rr}i e^{it\tau}d\mu_2(\tau)+y\int_{\rr}i\tau e^{it\tau}d\mu_2(\tau).$$
The last expression gives exactly the general solution as anticipated above. One just has to consider arbitrary Radon measures $d\mu_1(\tau)=\hat{A}(\tau)d\tau$ and
$d\mu_2(\tau)=\hat{B}(\tau)d\tau$ where $\hat{A}$ and $\hat{B}$ are the Fourier transforms of the two arbitrary functions $A$ and $B$.

}\end{example}

\end{document}